\documentclass[12pt,leqno]{article}
\usepackage{latexsym}

\newtheorem{defn}{Definition}[section]
\newtheorem{thm}[defn]{Theorem}
\newtheorem{cor}[defn]{Corollary}
\newtheorem{prop}[defn]{Proposition}
\newtheorem{lem}[defn]{Lemma}

\def\F{\mbox{${\cal F}$}}
\newcommand{\rr}{\mbox{$I \hspace{-2.6mm} R$}}
\newcommand{\cc}{\mbox{$I \hspace{-2.6mm} C$}}
\newcommand{\zz}{\mbox{$I \hspace{-2.6mm} Z$}}
\newcommand{\nn}{\mbox{$I \hspace{-2.6mm} N$}}
\input{arrow.tex}

\begin{document}
\title{\bf Uniqueness of the Fr\'{e}chet algebra topology on certain Fr\'{e}chet algebras}

\author{S. R. PATEL}

\maketitle
\vspace{1.5cm}
\begin{tabbing}
\vspace{1.5cm}
Department of Mathematics\\
Faculty of Engineering \& Technology\\
Marwadi Education Foundation's Group of Institutions\\
Rajkot, Gujarat, INDIA\\
Telephone (Res.):
+91 79 27418845\\
E-mails: srpatel.math@gmail.com, coolpatel1@yahoo.com
\> Secondary 13F25, 46H40\\
\end{tabbing}
\noindent {\bf Abstract.} In 1978, Dales posed a question about
the uniqueness of the $(F)$-algebra topology for $(F)$-algebras of
power series in $k$ indeterminates. We settle this in the
affirmative for Fr\'{e}chet algebras of power series in $k$
indeterminates. The proof goes via first completely characterizing
these algebras; in particular, it is shown that the
Beurling-Fr\'{e}chet algebras of semiweight type do not satisfy a
certain equicontinuity condition due to Loy. Some applications to
the theory of automatic continuity are also given, in particular
the case of Fr\'{e}chet algebras of power series in infinitely
many indeterminates.\\


\newpage
\section{Introduction.} Throughout the paper, ``algebra" will mean a
complex, commutative algebra with identity unless otherwise
specified. A {\it Fr\'{e}chet algebra} is a complete, metrizable
locally convex algebra $A$ whose topology $\tau$ may be defined by
an increasing sequence $(p_m)_{m\geq 1}$ of submultiplicative
seminorms. We may refer to $\tau$ as ``the Fr\'{e}chet topology of
$A$" in the following. The principal tool for studying Fr\'{e}chet
algebras is the Arens-Michael representation, in which $A$ is
given by an inverse limit of Banach algebras $A_m$ (see ~\cite[\S
5]{11} or ~\cite[\S 2]{12}).

Let $k\,\in\,\nn$ be fixed. We write $\F_k$ for the algebra
$\cc[[X_1,\,X_2,\,\dots ,\,X_k]]$ of all formal power series in
$k$ commuting indeterminates $X_1,\,X_2,\,\dots ,\,X_k$, with
complex coefficients. A fuller description of this algebra is
given in ~\cite[\S 1.6]{3}; we briefly recall some notation, which
will be used throughout the paper. Let $k\,\in\,\nn$, and let
$J\,=\,(j_1,\,j_2,\,\dots ,\,j_k)\,\in\,\zz^{+k}$. Set $$\mid
J\mid\,=\,j_1 + j_2 + \cdots + j_k;$$ ordering and addition in
$\zz^{+k}$ will always be component-wise. A generic element of
$\F_k$ is denoted by
$$\sum_{J\in\zz^{+k}}\lambda_J\,X^J\,=\,\sum
\{\lambda_{(j_1,\,j_2,\,\dots ,\,j_k)}\,X_1^{j_1}X_2^{j_2}\cdots
X_k^{j_k}\,:\,(j_1,\,j_2,\,\dots ,\,j_k)\,\in\,\zz^{+k}\}.$$ The
algebra $\F_k$ is a Fr\'{e}chet algebra when endowed with the weak
topology $\tau_c$ defined by the coordinate projections
$$\pi_I\,:\,\sum_{J\in \zz^{+k}}\lambda_J\,X^J\,\mapsto\,\lambda_I,\;\F_k\,\rightarrow \,\cc,$$
for each $I\,\in\,\zz^{+k}$. A defining sequence of seminorms for
$\F_k$ is $(p_m^{'})$, where $$p_m^{'}(\sum_{J\in
\zz^{+k}}\lambda_J\,X^J)\,=\,\sum_{\mid J\mid \leq m}\mid
\lambda_J \mid\;(m\,\in\,\nn).$$ A {\it Fr\'{e}chet algebra of
power series in $k$ variables} (shortly: FrAPS in $\F_k$) is a
subalgebra $A$ of $\F_k$ such that $A$ is a Fr\'{e}chet algebra
containing the indeterminates $X_1,\,X_2,\,\dots ,\,X_k$ and such
that the inclusion map $A\,\hookrightarrow \,\F_k$ is continuous
(equivalently, the projections $\pi_I,\,I\,\in\,\zz^{+k}$, are
continuous linear functionals on $A$). It is worthwhile mentioning
that in ~\cite[Corollaries 11.3 and 11.4]{4}, we show that the
time-honoured definitions of Banach and Fr\'{e}chet (and, more
generally, $(F)$-) algebras of power series in $\F_1$ contain a
redundant clause of the continuity of coordinate projections; this
is not possible in the several-variable case by ~\cite[Theorem
12.3]{4}.

Though Fr\'{e}chet algebras of power series in $k$ indeterminates
have been considered by Loy ~\cite{9}, recently these algebras --
and more generally, the power series ideas in general Fr\'{e}chet
algebras -- have acquired significance in understanding the
structure of Fr\'{e}chet algebras ~\cite{1, 3, 4, 12, 13, 14, 16}.
Thus it is of interest to investigate the following:
\begin{enumerate}
\item[{\rm (1)}] whether one can {\it completely} characterize
these algebras, \item[{\rm (2)}] whether such algebras have a
unique topology as Fr\'{e}chet algebras.
\end{enumerate}

In this paper we shall be concerned with the solution to the above
problems; our argument here is kept short because it uses key
ideas involved in the solution to these problems for the case
$k\,=\,1$ ~\cite{12}. (See Theorem 3.1 and Corollary 4.3 below.)
In Section 3, we obtain several results of independent interest.
Precisely, we shall classify FrAPS in $\F_k$ which do not satisfy
an equicontinuity condition (E): there is a sequence
$(\gamma_K)_{K\in\nn^k}$ of positive reals such that
$(\gamma_K^{-1}\pi_K)$ is equicontinuous ~\cite{9}. (See Theorem
3.10 below.)

We remark that the uniqueness of the Fr\'{e}chet topology of
$\F_k$ for each $k\,\in\,\nn$ is established in ~\cite{3}, and the
general case has been open since 1978 ~\cite[Question 11]{2}. We
use the structure of the closed ideals and their powers to
establish the uniqueness of the Fr\'{e}chet topology of FrAPS in
$\F_k$, this is not known for the larger algebra
$\F_\infty\,=\,\cc[[X_1,\,X_2,\,\dots ]]$ ~\cite{14} and FrAPS in
$\F_\infty$, and so we cannot apply our approach to establish the
uniqueness of the Fr\'{e}chet topology of FrAPS in $\F_\infty$.
Not only this, but Read also showed in ~\cite{14} that in the
absence of the uniqueness of the Fr\'{e}chet topology, the
Singer-Wermer conjecture cannot be established in the case of
Fr\'{e}chet algebras and thus, the situation on Fr\'{e}chet
algebras is markedly different from that on Banach algebras.
However, we shall give some remarks, establishing the uniqueness
of the Fr\'{e}chet topology of FrAPS in $\F_\infty$ admitting a
continuous norm. The solution of Question 11 of \cite{2} does
include the uniqueness of the Fr\'{e}chet topology of FrAPS in
$\F_1$, established in \cite{12}, as a special case.

A Fr\'{e}chet algebra $(A,\,(p_m))$ is said to be a {\it
Fr\'{e}chet algebra with power series generators}
$x_1,\,x_2,\,\dots ,\,x_k$ if each $y\,\in\,A$ is of the form
$$y\,=\,\sum_{J\in \zz^{+k}}\lambda_J\,x^J\,=\,\sum
\{\lambda_{(j_1,\,j_2,\,\dots ,\,j_k)}\,x_1^{j_1}x_2^{j_2}\cdots
x_k^{j_k}\,:\,(j_1,\,j_2,\,\dots ,\,j_k)\,\in\,\zz^{+k}\},$$ for
$\lambda_J$ complex scalars such that $\sum_{J\in \zz^{+k}}\mid
\lambda_J\mid \,p_m(x^J)\,<\,\infty$ for all $m$. Thus if $A$ is a
Fr\'{e}chet algebra with finitely many power series generators
$x_1,\,x_2,\,\dots,\,x_k$, then $A$ is a commutative, separable,
finitely generated Fr\'{e}chet algebra generated by
$x_1,\,x_2,\,\dots ,\,x_k$. A {\it semiweight function} on
$\zz^{+k}$ is a function $\omega\,:\,\zz^{+k}\,\rightarrow\,\rr$
such that $$\omega(M + N) \leq \omega(M)\omega(N), \;\omega(0) =
1\; \textrm{and}\; \omega(N) \geq 0\;(M, N \in \zz^{+k});$$ a
semiweight function is a {\it weight function} if for all $N \in
\zz^{+k}, \,\omega(N)
> 0$. To answer (2) above, we will investigate the following two
questions on FrAPS $A$ in $\F_k$ as an intermediate step.
\begin{enumerate}
\item[{\rm (A)}] When are $X_1,\,X_2,\,\dots ,\,X_k$ power series
generators for $A$?
\item[{\rm (B)}] When is $A$ isomorphic to an inverse limit of
Banach algebras of power series in $k$ variables?
\end{enumerate}

The solution to the above problems are given in Theorems
~\ref{Theorem 2.1_BP} and 3.10, respectively. In Section 4, we also
pose some interesting questions in automatic continuity theory.
\section{Fr\'{e}chet algebras.} Let $M$ be a closed maximal ideal of a Fr\'{e}chet algebra $A$. We
shall suppose from now on that
$\textrm{dim}(M/\overline{M^{2}})\,=\,k$ is finite (it is easy to
see that for finitely generated Fr\'{e}chet algebras this
condition is automatically satisfied; see ~\cite[Proposition
2.2]{15} for the Banach case). Then, by the remark following
Theorem 2.3 of ~\cite{15}, for each $n\,\in\,\nn$, the homogeneous
monomials of degree $n$ in $t_1,\,t_2,\,\dots ,\,t_k\,\in\,M$ are
representatives of a basis for $\overline{M^{n}}/\overline{M^{n +
1}}$ if and only if $\textrm{dim}(\overline{M^{n}}/\overline{M^{n
+ 1}})\,=\,C_{n+k-1,n}\,=\,\frac{(n+k-1)!}{n!(k-1)!}$ for all $n$,
and thus $M$ is not nilpotent. Thus, in a special case, we have
the following, with an eye on ~\cite[Lemma 2.1]{12}.
\begin{prop} \label{Proposition 2.3_SM} Let $(A,\,(p_m))$ be a
commutative, unital Fr\'{e}chet algebra with the Arens-Michael
isomorphism $A\,\cong
\,\lim\limits_{\longleftarrow}(A_{m};\,d_{m})$. Suppose that there
exists a fixed $k\,\in\,\nn$ such that $M$ is a closed maximal
ideal of $A$ such that: (i) $\bigcap _{n\geq 1}\overline{M^{n}} =
\{0\}$ and (ii) $\textrm{dim}(\overline{M^{n}}/\overline{M^{n +
1}}) = C_{n+k-1,n}$ for all $n$. Then there exist
$t_1,\,t_2,\,\dots ,\,t_k\,\in\,M$ such that
$\overline{M^{n}}\,=\,\overline{M^{n + 1}}\;\oplus
\,\textrm{sp}\{t^I\,:\,\mid I\mid \,=\,n\}$ for each $n\,\geq\,1$.
Assume further that each $p_m$ is a norm. Then, for each
sufficiently large $m$, $M_m$ is a non-nilpotent maximal ideal of
$A_m$ such that: (a) $\bigcap _{n\geq 1}\overline{M_{m}^{n}} =
\{0\}$ and (b) $\textrm{dim}(\overline{M_m^{n}}/\overline{M_m^{n +
1}})\,=\,C_{n+k-1,n}$ for all $n$.
\end{prop}
{\it Proof.} The first half of the proof has already been
discussed above. For the second half of the proof, follow
~\cite[Proposition 2.3]{13}.

Concerning Proposition 2.1, the counter-examples (see ~\cite{13};
for the one-variable case) show that the assumption that each
$p_m$ is a norm on $A$ cannot be dropped. The algebra $\F_k$ is a
trivial counter-example in the several-variable case. We also
remark that, in the case where
$\textrm{dim}(M/\overline{M^{2}})\,=\,1$, one deduces
$\textrm{dim}(\overline{M^{n}}/\overline{M^{n + 1}})\,=\,1$ for
all $n$ in ~\cite[Proposition 2.3]{12}, and so we do not require
$\textrm{dim}(\overline{M^{n}}/\overline{M^{n + 1}})\,=\,1$ for
all $n$ as a stronger hypothesis, but then we do require $M$ to be
non-nilpotent there. Below, we exhibit an easy counter-example to
show that the hypothesis that
$\textrm{dim}(\overline{M^{n}}/\overline{M^{n +
1}})\,=\,C_{n+k-1,n}$ for all $n\,\in\,\nn$ is not redundant in
the proposition above (many thanks to Professor H. G. Dales for
calling my attention to this counter-example).

Let $$B\, =\, \F_2 =\, \cc[[X,\,Y]],$$ with the usual Fr\'{e}chet
algebra topology $\tau_c$ and let $J$ be the ideal generated by
the element $X^2 - Y^3$. Since $B$ is noetherian ~\cite[VII,
Corollary p. 139 and Theorem 4{'}]{17}, all ideals in $B$ are
closed by ~\cite[Theorem 5]{19}, so $J$ is closed. Hence the
quotient $A = B/J$ is a noetherian Fr\'{e}chet algebra, with all
the ideals in $A$ closed. Clearly, $M/J$ is the unique maximal
ideal in $A$, where $M\,=\,\ker \pi_0\;(0\,=\,\{0,\,0\})$, in $B$.
The two elements $X + J$ and $Y + J$ are linearly independent
modulo $(M/J)^2$, and so dim$((M/J)/(M/J)^2)\,=\,2$ since
$M/M^2\,\cong \,(M/J)/(M/J)^2$. However $(X + J)^2 \,\in
\,(M/J)^3$, so dim$((M/J)^2/(M/J)^3)\,=\,2$ since $XY + J$ and
$Y^2 + J$ are linearly independent modulo $(M/J)^3$.

Next, to see that this is a counter-example, we show that $\bigcap
_{n\geq 1}(M/J)^n \,\neq \, J$, the zero element of $A$. To see
this, let us start with an element $g$ of $B$ such that
$$g\,\in\,p X^2 + q XY + r Y^2 + M^3 + J,$$ for some $p, \,q,\, r\, \in\,
\cc$, that is, $$g\,\in\,\cc X^2 + \cc XY + \cc Y^2 + M^3 + J$$
and suppose that $g\,\in\,M^4 + J$. Then
$$g\,\in\,(X^2 - Y^3)(a + b X +c Y) + M^4,$$ for some $a, b, c \in \cc$
because all other terms in $J\,=\,(X^2 - Y^3)B$ are in $M^4$. Thus
we can see that there exist $a_1, b_1, c_1 \in \cc$ such that
$$p X^2 + q XY + r Y^2 = (X^2 - Y^3)(a_1 + b_1 X + c_1 Y) + M^3.$$ Now,
equating the coefficients of $XY$ and $Y^2$, we see that $q = r =
0$, and, equating the coefficients of $X^2$, we see that $p =
a_1$, and then, equate the coefficients of $Y^3$ to see that $0 =
a_1$. Thus $p = 0$. We conclude that $M^3 + J\,=\,M^4 + J$. One
can generalize this idea to see that $M^n + J\,=\,M^{n+1} + J$ for
each $n \geq 3$, the only element of $B$ that belongs to $M^n + J$
is actually in $M^{n +1} + J$. Thus we conclude that $M^3 +
J\,=\,M^n + J $ for each $n \geq 3$, and so $\bigcap _{n\geq
1}(M/J)^n = (M/J)^3\,\neq\,J$.
\section{Fr\'{e}chet algebras of power series in $\F_k$.} We now
turn to the problem of describing {\it all} those commutative
Fr\'{e}chet algebras which may be continuously embedded in $\F_k$
in such a way that they contain the polynomials in
$X_1,\,X_2,\,\dots ,\,X_k$. The following theorem completely
characterizes separable FrAPS in $\F_k$. The method of proof will
be used again in the proof of Theorem 3.10.
\begin{thm} \label{Theorem 3.1_SM} Let $A$ be a commutative, unital Fr\'{e}chet
algebra. Suppose that there exists a fixed $k\,\in\,\nn$ such that
$A$ contains a closed maximal ideal $M$ such that: (i) $\bigcap
_{n\geq 1}\overline{M^{n}} =$ $\{0\}$; and (ii)
dim$(\overline{M^n}/\overline{M^{n+1}}) = C_{n+k-1,n}$ for all
$n$. Then $A$ is a Fr\'{e}chet algebra of power series in $\F_k$.
The converse holds if the polynomials in $X_1,X_2,\dots ,X_k$ are
dense in $A$.
\end{thm}
{\it Proof.} The proof is similar to that of ~\cite[Theorem
3.1]{12}, and will be \linebreak outlined only. Supposing $A$
satisfies the stated conditions, there exist \linebreak $t_1, t_2,
\dots , t_k$ $\in M$ such that
$$\overline{M^{n}}\,=\,\overline{M^{n + 1}}\;\oplus \,\textrm{sp}\{t^I\,:\,\mid I\mid \,=\,n\},$$
for each $n\,\geq\,1$, by Proposition ~\ref{Proposition 2.3_SM}.
Let $x\,\in\,A$. Then a simple induction on $n$ shows that for
$n\,\geq\,1$, $ x\,=\,\sum_{\mid I\mid\,\leq\,n}\lambda
_{I}\,t^{I} \,+\, y_{n},$ where $y_{n}\,\in\,\overline{M^{n + 1}}$
and the $(\lambda_{I})$ are uniquely determined. Hence the
functionals $\pi_{J}\,:\,x\,\mapsto \,\lambda_{J}$ are uniquely
defined, and linear for all $J\,\in\,\nn^k$. If
$x\,\in\,\textrm{ker}\,\pi_{J}$ for all $J\,\in\,\nn^k$, then
$x\,\in\,\bigcap _{n\geq 1}\overline{M^{n}}\,=\,\{0\}$. Thus the
mapping $$x\,\mapsto \,\sum_{I\,\in\,\zz^{+k}}\pi_{I}(x)\,t^{I}$$
is an isomorphism of $A$ onto an algebra of formal power series in
$\F_k$.

\indent Carrying over the topology via this isomorphism, the
result will follow once we show that the functionals ${\pi_{J}}$
are continuous for each $J$. Clearly $\pi_{0,\dots,\,0}$ is
continuous since $M\,=\,\ker \pi_{0,\dots,\,0}$ is a closed
maximal ideal of $A$. Let $J\,\in\,\nn^k$, and assume that
$\pi_{I}$ is continuous for each $\mid I \mid \,<\,\mid J \mid$,
and take $(x_{n})$ in $A$ with $x_{n}\,\rightarrow\,0$ as
$n\,\rightarrow\,\infty$. Then
$$x_{n}\,=\,\sum_{\mid I\mid\,\leq\,\mid J\mid}\pi_{I}(x_{n})\,t^{I}\,+\,y_{n,\, \mid J\mid},$$
for some $y_{n,\,\mid J\mid}\,\in\,\overline{M^{\mid J\mid + 1}}$.
It follows that $$\sum_{\mid I\mid\,=\,\mid
J\mid}\pi_I(x_{n})\,t^{I}\,+\,y_{n,\, \mid
J\mid}\,\rightarrow\,0,$$ so if $\pi_I(x_{n})$ does not converge
to $0$ for at least one $I$ ($\mid I\mid\,=\,\mid J\mid$) we
deduce that some non-zero linear combination of $t^I,\,\mid
I\mid\,=\,\mid J\mid$, lies in $\overline{M^{\mid J\mid + 1}}$, a
contradiction. Thus each $\pi_J$ is continuous.

Conversely, let $A$ be a Fr\'{e}chet algebra of power series in
$\F_k$ such that the polynomials in $X_1,X_2,\dots ,X_k$ are dense
in $A$. Setting $M\,=\,\textrm{ker}\,\pi_{0,\dots,\,0}$, we have
$\overline{M^n}\, \subset \, \textrm{ker}\,\pi_N$ for each $N\,\in
\zz^{+k}$ with $\mid N\mid\,=\,n - 1\;(n\,\in\,\nn)$, so that
$\bigcap _{n\geq 1}\overline{M^{n}}\, =\, \{0\}$. Clearly
$$\overline{M^{n+1}}\, \oplus \,\textrm{sp}\{t^I\,:\,\mid I\mid \,=\,n\}\,\subseteq \,\overline{M^n}.$$
Let $$M_{n+1}\, =\, \{a \,\in\,A\,:\,o(a)\,\geq\,n + 1\}.$$ Since
$M_{n+1}$ is closed by continuity of the functionals $\pi_J$, we
have \linebreak $\overline{M^{n+1}} \,\subseteq\, M_{n+1}$. Given
that the polynomials are dense in $A$, the series actually
converge in $A$, and so $\overline{M^{n+1}} \,=\, M_{n+1}$. Hence
dim$(\overline{M^n}/\overline{M^{n+1}})\,=\,C_{n+k-1,n}$ for all
$n$, and the theorem follows.

We note that, in {\it all} FrAPS $A$ in $\F_k$, we have
$M\,=\,\ker \pi_{0,\dots,\,0}$ is a non-nilpotent, closed maximal
ideal such that $\bigcap _{n\geq 1}\overline{M^{n}}\, =\, \{0\}$.
Two counter-examples in the one-variable case (see ~\cite[Remarks
1 (b)]{12}) show that the assumption that the polynomials are
dense in $A$ cannot be dropped in the above theorem.

We now turn to answer (A) above. The following lemma, whose proof
we omit, is a several-variable-analogue of ~\cite[Lemma 3.2]{12}
(for proof, see ~\cite[Lemma 2.2]{1}. Define a subset $A_1$ of $A$
by
$$A_1\,:=\,\{y\,\in\,A\,:\,\sum_{J\,\in\,\zz^{+k}}\mid \lambda_J\mid
p_m(X^J)\,<\,\infty\;\textrm{for all}\; m\}$$ in the lemma.
\begin{lem} \label{Lemma 3.2_SM} Let $A$ be a Fr\'{e}chet algebra of power series in $\F_k$.
Then:
\begin{enumerate}
\item[{\rm (1)}] $A_1$ is continuously embedded in $A$; \item[{\rm
(2)}] $A_1$ is a Fr\'{e}chet algebra having power series generators
$X_1,\,X_2,\,\dots ,\,X_k$; \item[{\rm (3)}] $A_1$ is a Banach
algebra provided that $A$ is a Banach algebra.
\end{enumerate}
\end{lem}

We recall that elements $x_1,\,x_2,\,\dots ,\,x_k$ in a
Fr\'{e}chet algebra $A$ generate a {\it multi-cyclic basis} if
each $y\,\in\,A$ can be uniquely expressed as
$$y\,=\,\sum_{J\in \zz^{+k}}\lambda_J\,x^J\,=\,\sum
\{\lambda_{(j_1,\,j_2,\,\dots ,\,j_k)}\,x_1^{j_1}x_2^{j_2}\cdots
x_k^{j_k}\,:\,(j_1,\,j_2,\,\dots ,\,j_k)\,\in\,\zz^{+k}\},$$
$\lambda_J$ complex scalars. A seminorm $p$ on a Fr\'{e}chet
algebra $A$, having power series generators $x_1,\,x_2,\,\dots
,\,x_k$ generating a multi-cyclic basis for $A$, is a {\it power
series seminorm} if $$p(y\,=\,\sum_{J\in
\zz^{+k}}\lambda_J\,x^J)\,=\,\sum_{J\,\in\,\zz^{+k}}\mid
\lambda_J\mid p(x^J)\;(y\,\in\,A).$$
\begin{cor} \label{Corollary 3.3_BP} Let $A$ be a Fr\'{e}chet algebra $A$
having power series generators $x_1,\,x_2,\,\dots ,\,x_k$. Then
$x_1,\,x_2,\,\dots ,\,x_k$ generate a multi-cyclic basis for $A$
if and only if the topology of $A$ is defined by a sequence of
power series seminorms.
\end{cor}
{\it Proof.} For the proof of ``only if" part, follow ~\cite[Lemma
2.2]{1}.

Next, we define {\it Beurling-Fr\'{e}chet algebras
$\ell^1(\zz^{+k},\,\Omega)$ of semiweight type}, and list some of
their useful properties.

First, we recall from Introduction that $\omega$ is a {\it proper
semiweight} if \linebreak $\omega(N_0)\,=\,0$ for some
$N_{0}\,\in\, \nn^k$. Let $k\,\in\,\nn$, and let $(A, \,(p_m))$ be
the Beurling-Fr\'{e}chet algebra
$$\ell^1(\zz^{+k},\,\Omega) := \{f=\sum_{J\in \zz^{+k}}\lambda_J
X^J\in \F_k:\sum_{J\in \zz^{+k}}\mid\lambda_J\mid\
\omega_m(J)<\infty\,\textrm{for all}\,m\},$$ where
$\Omega\,=\,(\omega_m)$ is a separating and increasing sequence of
semiweight functions on $\zz^{+k}$ defined by
$\omega_m(J)\,=\,p_m(X^J)$. If $\Omega$ is an increasing sequence
of weight functions on $\zz^{+k}$, then we define
$$\rho=\sup_m\rho_m, \;\textrm{where}\; \rho_m=\inf_{N\in
\zz^{+k}}\omega_m(N)^{\frac{1}{\mid N\mid}}.$$ Thus, $\rho\,=\,0$
if and only if $\rho_m\,=\,0$ for each $m$, if and only if for
each $m$ $\ell^1(\zz^{+k},\,\omega_m)$ is a local Banach algebra
in the Arens-Michael representation of $\ell^1(\zz^{+k},\,\Omega)$
if and only if $\ell^1(\zz^{+k},\,\Omega)$ is a local Fr\'{e}chet
algebra, and $\rho\,>\,0$ if and only if $\rho_m\,>\,0$ for some
$m$ if and only if for each $l\,\geq\,m$,
$\ell^1(\zz^{+k},\,\omega_l)$ is a semisimple Banach algebra in
the Arens-Michael representation of $\ell^1(\zz^{+k},\,\Omega)$ if
and only if $\ell^1(\zz^{+k},\,\Omega)$ is a semisimple
Fr\'{e}chet algebra.

Suppose that $\Omega$ is a separating and increasing sequence of
{\it proper semiweights} on $\zz^{+k}$. Then $\rho \,= 0$ if and
only if $\ell^1(\zz^{+k},\,\Omega)$ is a local Fr\'{e}chet algebra
\linebreak if and only if the completion of
$\ell^1(\zz^{+k},\,\omega_m)/\ker p_m$ under the induced norm
$p_m$ is a local Banach algebra for all $m$. In this case,
$\ell^1(\zz^{+k},\,\Omega)$ is either $\F_k$ or a local FrAPS in
$\F_k$. We call such a Beurling-Fr\'{e}chet algebra
$\ell^1(\zz^{+k},\,\Omega)$ an {\it algebra of semiweight type}.
We note that the unique maximal ideal of
$\ell^1(\zz^{+k},\,\Omega)$ is
$$\{f=\sum_{J\in \zz^{+k}}\lambda_J\,X^J\in
\ell^1(\zz^{+k},\,\Omega)\,:\,\lambda_0\,=\,0\}.$$ For example, if
$k = 1$, then, by ~\cite[Theorem 2.1]{1},
$\ell^1(\zz^+,\,\Omega)\,=\,\F$, which is an inverse-limit of
finite-dimensional algebras, and is also a local algebra. In this
case, $$\omega_m\,:\,\zz^+\,\rightarrow
\,[0,\,\infty),\;\omega_m(n)\,=\,p_m(X^n),$$ is $1$, if $n\leq m$
and is $0$, if $n > m$ (~\cite[p. 131]{4}). If $k = 2$, then, by
Theorem 3.4 \linebreak below, $\ell^1(\zz^{+2},\,\Omega)$ is
either $\F_2$ or $A_X$ or $A_Y$ (all the three algebras are
local), where
$$A_X\,:=\,\{f\,=\,\sum_{i,j}\lambda_{i,j}X^iY^j\,\in\,\F_2\,:\,p_m(f)\,:=\,
\sum_{j=0}^{\infty}\sum_{i=0}^m\,\mid\lambda_{i,j}\mid\,<\,\infty\;\textrm{for
all}\;m\}$$ in which case
$$\omega_m\,:\,\zz^{+2}\,\rightarrow\,[0,\,\infty),\;\omega_m(i,j)\,=\,p_m(X^iY^j),$$
is $1$, if $i\,\leq\,m,\,j\,\in\,\zz^+$ and is $0$, if
$i\,>\,m,\,j\,\in\,\zz^+$, and where
$$A_Y\,:=\,\{f\,=\,\sum_{i,j}\lambda_{i,j}X^iY^j\,\in\,\F_2\,:\,p_m(f)\,:=\,
\sum_{i=0}^{\infty}\sum_{j=0}^m\,\mid\lambda_{i,j}\mid\,<\,\infty\;\textrm{for
all}\;m\}$$ in which case
$$\omega_m\,:\,\zz^{+2}\,\rightarrow\,[0,\,\infty),\;\omega_m(i,j)\,=\,p_m(X^iY^j),$$
is $1$, if $j\,\leq\,m,\,i\,\in\,\zz^+$ and is $0$, if
$j\,>\,m,\,i\,\in\,\zz^+$. For
$\ell^1(\zz^{+k},\,\Omega)\,=\,\F_2$, with
$$p_m(f)\,:=\,\sum_{0\leq i+j \leq m} \mid \lambda_{i,j}\mid,$$ we
define $\Omega\,=\,(\omega_m)$, where
$$\omega_m\,:\,\zz^{+2}\,\rightarrow\,[0,\,\infty),\;\omega_m(i,j)\,=\,p_m(X^iY^j),$$
is $1$, if $0\,\leq\,i+j\,\leq\,m$ and is $0$, if $i+j\,>\,m$
(~\cite[p. 131]{4}). Clearly, $A_X \cong A_Y$ under the
interchange of variables $X$ and $Y$.

If $k\,=\,3$, then, again by Theorem 3.4 below,
$\ell^1(\zz^{+3},\,\Omega)$ is either $\F_3$ or $A_X$ or $A_Y$ or
$A_Z$ or $A_{X,Y}$ or $A_{X,Z}$ or $A_{Y,Z}$ defined analogously
(in fact, $A_X\,\cong\,A_Y\,\cong\,A_Z$ and
$A_{X,Y}\,\cong\,A_{X,Z}\,\cong\,A_{Y,Z}$). We can extend above
arguments for $k\,\geq\,4$. Thus, for $k\,\in\,\zz^+$, we have
Beurling-Fr\'{e}chet algebras $\ell^1(\zz^{+k},\,\Omega)$ of
semiweight type, with the following properties:

(1) $\F_k$ is the {\it only} Fr\'{e}chet algebra of finite type
among FrAPS in $\F_k$, by Corollary 3.8 below; rest of
$\ell^1(\zz^{+k},\,\Omega)$ are not Fr\'{e}chet algebras of finite
type (see ~\cite{7}).

(2) The Arens-Michael representations of
$\ell^1(\zz^{+k},\,\Omega)$ do not contain BAPS in $\F_k$ (for the
time being, we assume that such algebras have unique Fr\'{e}chet
topology, which we shall prove later); for proof, see Remark
preceding to Corollary \ref{Corollary 3.2_BP} below.

(3) The polynomials in $k$ variables are dense in
$\ell^1(\zz^{+k},\,\Omega)$.

(4) The Fr\'{e}chet topology $\tau$ of $\ell^1(\zz^{+k},\,\Omega)
\;(\neq\,\F_k)$, defined by a sequence $(p_m)$, is finer than
$\tau_c$ of $\F_k$, but surely not equivalent otherwise the
$\tau$-closure of the algebra of polynomials in $k$ variables
(which is $\ell^1(\zz^{+k},\,\Omega)$) is equal to $\F_k$, a
contradiction. One can also deduce a contradiction from the
statement (1) above. Hence the rest of $\ell^1(\zz^{+k},\,\Omega)$
differ from $\F_k$ (also, from the statement (1) point of view as
well).

(5) $\F_r,\,1\leq r\leq k - 1$, can be regarded as closed
subalgebras of $\ell^1(\zz^{+k},\,\Omega)$ via the obvious
quotient maps (e.g., if $k\,=\,2$, then $\cc[[X]]\,=\,\F_1$ can be
regarded as a closed subalgebra of $A_X$; the quotient map from
$A_X$ obtained by setting $Y\,=\,0$ is denoted by
$$\pi\,:\,\sum_{i,j}\lambda_{i,j}X^iY^j\,\mapsto\,\sum_{i=0}^{\infty}\lambda_{i,0}X^i,\;A_X\,\rightarrow\,\F_1).$$
Hence all Beurling-Fr\'{e}chet algebras
$\ell^1(\zz^{+k},\,\Omega)$ of semiweight type are local
Fr\'{e}chet algebras since the closed subalgebras $\F_r,\,1\leq
r\leq k,$ are local Fr\'{e}chet algebras, and the unique maximal
ideal $M$ is
$$\{f=\sum_{J\in \zz^{+k}}\lambda_J\,X^J\in
\ell^1(\zz^{+k},\,\Omega)\,:\,\lambda_0\,=\,0\}.$$ Also, for a
fixed $k\,\in\,\nn$, there are finitely many Beurling-Fr\'{e}chet
algebras $\ell^1(\zz^{+k},\,\Omega)$ of semiweight type, and these
algebras can be properly nested (for example, if $k\,=\,3$, then
$A_X\,\subset\,A_{X,Y}\,\subset\,\F_3$). Further, if $(A,\,(q_m))$
is a FrAPS in $\F_k$ such that the $q_m$ are proper seminorms on
$A$, then $A$ is continuously embedded in the ``least"
Beurling-Fr\'{e}chet algebra $\ell^1(\zz^{+k},\,\Omega)$ of
semiweight type (note that there might be several such
Beurling-Fr\'{e}chet algebras of semiweight type containing $A$).
Moreover it is clear that if such $A$ contains a
Beurling-Fr\'{e}chet algebra of semiweight type such that
$$\ell^1(\zz^{+k},\,\Omega_1)\,\hookrightarrow\,A\,\hookrightarrow\,\ell^1(\zz^{+k},\,\Omega_2)$$
continuously, then, depending on the $q_m$, $A$ is either
$\ell^1(\zz^{+k},\Omega_1)$ or $\ell^1(\zz^{+k},\Omega_2)$ or none
of these in which case both the inclusions are proper (for
example, $A_X\,\hookrightarrow\,A_{X,Y}\,\hookrightarrow\,\F_3$).
We use these facts in the proof of Theorem 3.7 below.

We call $\ell^1(\zz^{+k},\,\Omega)$ a {\it Beurling-Fr\'{e}chet
algebra of weight type} if $\Omega$ is an increasing sequence of
weight functions $\omega_m$ on $\zz^{+k}$. In this case, the
topology $\tau$ is defined by the increasing sequence $(p_m)$ of
norms defined in terms of the $\omega_m$, and the corresponding
Arens-Michael representation contains Beurling-Banach algebras
$\ell^1(\zz^{+k},\,\omega _m)$. When it can cause no confusion, we
may call $\ell^1(\zz^{+k},\,\Omega)$ a Beurling-Fr\'{e}chet
algebra of (semi)weight type if $\Omega$ is an increasing sequence
of semiweight functions $\omega_m$ on $\zz^{+k}$ (which would
include both type of algebras: algebras of semiweight type as well
as algebras of weight type).
\begin{thm} \label{Theorem 2.1_BP} Let $A$ be a Fr\'{e}chet algebra
 of power series in $\F_k$. Suppose that
$X_1,\,X_2,\,\dots ,\,X_k$ are power series generators for $A$.
Then $A$ is the Beurling-Fr\'{e}chet algebra
$\ell^1(\zz^{+k},\,\Omega)$ for an increasing sequence $\Omega$ of
semiweight functions on $\zz^{+k}$.
\end{thm}
{\it Proof.} Let $A$ satisfy the stated conditions. By the
uniqueness of the formal power series expression and the fact that
$X_1,\,X_2,\,\dots ,\,X_k$ are power series generators for $A$, it
follows that $\{X^J\,:\,J\,\in\,\zz^{+k}\}$ is a multi-cyclic
basis for $A$. By Corollary \ref{Corollary 3.3_BP}, the
Fr\'{e}chet topology $\tau$ of $A$ is defined by an increasing
sequence $(p_m)$ of power series seminorms; and for each
$y\,=\,\sum_{J\,\in\,\zz^{+k}}\lambda_JX^J$ in $A$, we have
$y\,=\,\lim_n\,\sum_{\mid J\mid\,\leq\,n}\lambda_JX^J$ in the
topology $\tau$. For each $m$,
$$\omega_m\,:\,\zz^{+k}\,\rightarrow\,[0,\,\infty),\;\omega_m(N)\,=\,p_m(X^N)$$
define a separating sequence $\Omega$ of semiweight functions. Let
$\ell^1(\zz^{+k},\,\Omega)$ be as defined following Corollary 3.3.
Since each $p_m$ is a power series seminorm,
$A\,\subset\,\ell^1(\zz^{+k},\,\Omega)$. In fact,
$A\,=\,\ell^1(\zz^{+k},\,\Omega)$. For let
$$f\,=\,\sum_{J\,\in\,\zz^{+k}}\lambda_JX^J\,\in\,\ell^1(\zz^{+k},\,\Omega).$$
Let $f_n\,=\,\sum_{\mid J\mid\,\leq\,n}\lambda_JX^J$. Since
$X_i\,\in\,A$ for $i\,=\,1,\,2,\,\dots,\,k$, each $f_n\,\in\,A$; and
$(f_n)$ is a Cauchy sequence in $A$. Thus $f\,\in\,A$.

Now we have the following possibilities:
\begin{enumerate}
\item[{\rm (a)}] all $\omega_m$ are weights; \item[{\rm (b)}] no
$\omega_m$ is a weight; \item[{\rm (c)}] at least one
$\omega_{m_0}$ fails to be a weight.
\end{enumerate}
In the case (c), let $G\,=\,\{\omega_m\,:\,\omega_m\;\textrm{is
not a weight}\}$. If $G$ is finite the corresponding $p_m$ may be
deleted, if $G$ is infinite the corresponding $p_m$ can be taken
to define the topology so reducing consideration to the case (b).
Assume (b). Then for each $m$, there is $N \in \zz^{+k}$ such that
$\omega_m(N) = 0$. Then $$\omega_m(N + M) \leq \omega(N)\omega(M)
= 0\;(M \in \zz^{+k}).$$ Now, depending on $\Omega = (\omega_m)$,
we have $A = \ell^1(\zz^{+k},\,\Omega)$ a Beurling-Fr\'{e}chet
algebra of semiweight type. In the case (a), we have
$A\,=\,\ell^1(\zz^{+k},\,\Omega)$ a Beurling-Fr\'{e}chet algebra
of weight type. The theorem follows.

Next, let $A$ be a FrAPS in $\F_k$. We call a seminorm $p$ on $A$
{\it closable} if for any $p$-Cauchy sequence $(f_l)$ in $A$,
$f_l\,\rightarrow\,0$ in $\tau_c$ implies that
$p(f_l)\,\rightarrow\,0$. We define $p$ to be of type (E) if given
$M\,\in\,\zz^{+k}$, there exists $c_M\,>\,0$ such that $$\mid
\pi_M(f)\mid\, \leq \,c_Mp(f),$$ for all $f \,\in \,A$ ~\cite{9}.
A seminorm of type (E) is a norm. Also, closability of a norm on a
normed algebra of power series in $k$ indeterminates is a
necessary and sufficient condition for the completion to be a BAPS
in $\F_k$ (see ~\cite[Lemma 3.5]{1}.

We now answer (B) stated in the introduction. The following
proposition, whose proof we omit, is a several-variable-analogue of
~\cite[Proposition 3.1]{1}.
\begin{prop} \label{Proposition 3.1_BP} Let $A$ be a Fr\'{e}chet algebra
of power series in $\F_k$. Let $p$ be a continuous
submultiplicative seminorm on $A$. Let $\ker p\,=\,\{f\,\in\,
A\,:\,p(f)\,=\,0\}$. Let $A_p$ be the completion of $A/\ker p$ in
the norm $\|f + \ker p\|_p\,=\,p(f)$. Then the following are
equivalent.
\begin{enumerate}
\item[{\rm (i)}] $p$ is a norm and $A_p$ is a Banach algebra of
power series in $\F_k$. \item[{\rm (ii)}] $p$ is closable and of
type (E).
\end{enumerate}
\end{prop}
{\bf Remark.} The above proposition can be used to prove the
property (2) as follows: consider $A_Y\,\subset\,\F_2$. The
topology $\tau$ of $A_Y$ is given by a sequence $(p_m)$ of proper
power series seminorms, hence they are not of type (E). So,
$(A_Y)_m$ cannot be a BAPS in $\F_2$.
\begin{cor} \label{Corollary 3.2_BP} Let
$A\,=\,\lim\limits_{\longleftarrow}A_{m}$ be the Arens-Michael
representation of a Fr\'{e}chet algebra of power series in $\F_k$.
Assume that each $p_m$ is a norm. Then each $A_m$ is a Banach
algebra of power series in $\F_k$ if and only if each $p_m$ is a
closable norm of type (E).
\end{cor}

It is readily seen that a Fr\'{e}chet algebra of power series in
$\F_k$ satisfies Loy's condition (E) in ~\cite{9} if and only if
$A$ admits a continuous norm of type (E) \linebreak if and only if
the topology of $A$ is defined by a sequence of norms of type (E).
\linebreak \indent Next, we give characterizations of a
Beurling-Fr\'{e}chet algebra $\ell^1(\zz^{+k},\,\Omega)$ of
semiweight type. We have the following elementary, but crucial,
theorem. By identifying the series expansion in $x_1,\,x_2,\,\dots
,\,x_k$  with the series expansion in $X_1,\,X_2,\,\dots ,\,X_k$,
Fr\'{e}chet algebras with a multi-cyclic basis are realized as
Fr\'{e}chet algebras of power series in $\F_k$, the projections
$\pi_J$ being continuous. Note that by a {\it proper} seminorm we
mean a seminorm that is not a norm.
\begin{thm} \label{Theorem 3.3_SM} Let $A$ be a Fr\'{e}chet algebra
of power series in $\F_k$. Then $A$ is either a
Beurling-Fr\'{e}chet algebra $\ell^1(\zz^{+k},\,\Omega)$ of
semiweight type or the Fr\'{e}chet topology $\tau$ of $A$ is
defined by a sequence $(p_m)$ of norms.
\end{thm}
{\it Proof.} If $A$ is Banach, then certainly the topology $\tau$
of $A$ is defined by a norm, and so $A$ is not equal to a
Beurling-Fr\'{e}chet algebra $\ell^1(\zz^{+k},\,\Omega)$ of
semiweight type. Now suppose that $A$ is a non-Banach FrAPS in
$\F_k$. Let $(p_m)$ be an increasing sequence of seminorms
defining the Fr\'{e}chet topology $\tau$ of $A$, and set
$$G\,=\,\{l\,\in\,\nn\,:\,p_l\;\textrm{is a proper seminorm
on}\,A\}.$$ If $G$ is finite the corresponding $p_l$ may be
deleted and we have a new sequence of norms, defining the same
Fr\'{e}chet topology $\tau$ of $A$. Otherwise, $G$ is infinite the
corresponding $p_l$ can be taken to define the Fr\'{e}chet
topology $\tau$ of $A$. Then, by Lemma \ref{Lemma 3.2_SM}, there
exists a Fr\'{e}chet subalgebra $(A_1,\,(q_m))$ of $(A,\,(p_m))$
continuously embedded in $A$; $A_1$ is a Fr\'{e}chet algebra with
power series generators $X_1,\,X_2,\,\dots ,\,X_k$. By Theorem
\ref{Theorem 2.1_BP}, $A_1$ is a Beurling-Fr\'{e}chet algebra
$\ell^1(\zz^{+k},\,\Omega)$ of (semi)weight type. To rule out a
possibility of $A_1$ being a Beurling-Fr\'{e}chet algebra
$\ell^1(\zz^{+k},\,\Omega)$ of weight type, fix $m\,\geq\,1$. Then
there is a non-zero element $f\,\in\,A$ that belongs to $\ker
p_m$. Let $f_l$ be the {\it initial form} of $f$ ~\cite[p.
130]{17}. Thus $f_l$ is a homogeneous polynomial of degree $n(m)$
(which is minimal as $f$ varies in $\ker p_m$). Since $\ker p_m$
is a closed ideal in $A$ and the inclusion map
$A\,\hookrightarrow\,\F_k$ is continuous with $A$ a dense
subalgebra of $(\F_k,\,\tau_c)$, by ~\cite[Lemma B.10]{11}, the
closure of $\ker p_m$ (say $I$) in the topology $\tau_c$ is a
closed ideal in $\F_k$ such that $\ker p_m\,=\,I\,\bigcap\,A$ and
that $\ker p_m \,\subseteq \,I\,\subseteq \,\ker p_{l(m)}^{'}$.
Now, by the lemma from ~\cite[p. 136]{17}, there exists an
automorphism $\psi$ of $\F_k$ such that $\psi(f)\,=\,f_l$. Hence
$f_l\,\in\,\ker p_{l(m)}^{'}$ (note that this is true since we
have replaced $\F_k$ by $\psi(\F_k)\,=\,\F_k$, which is a FrAPS).
Since $\ker p_{l(m)}^{'}$ is finitely generated by the monomials
and $f_l$ is a homogeneous polynomial of degree $n(m)$ in $\ker
p_{l(m)}^{'}$, so it is, indeed, finitely generated by the
monomials of degree $n(m)$. Since $I$ is a finitely generated
ideal in $\F_k$ ($\F_k$ being noetherian) and $I\,\subseteq \,\ker
p_{l(m)}^{'}$, $I$ is finitely generated by the monomials of
degree $n(m)$. So $f_l\,\in\,I$. Evidently, $f_l\,\in\,A$, being a
homogeneous polynomial of degree $n(m)$. So $f_l\,\in\,\ker p_m$.

Hence
$$p_m(X^{n(m)})\,=q_m(X^{n(m)})\,=\,0\; (m\,\in\,\nn),$$ and so,
$A_1$ is, indeed, a Beurling-Fr\'{e}chet algebra
$\ell^1(\zz^{+k},\,\Omega)$ of semiweight type. Hence $A_1$ is a
local algebra, and therefore $A$ is also a local algebra. Now
there are several cases; we consider these cases for $k\,=\,2$
(one can modify the following arguments for any $k$). Thus $f_l$
is a monomial of degree $n(m)$ in $X$ and $Y$. We have the
following three cases.

Case 1. For each $m$, $\ker p_m$ is finitely generated by the
monomials of \linebreak degree $n(m)$ in $X$ and $Y$ in which case
$I$ is also finitely generated by the \linebreak monomials of
degree $n(m)$ in $X$ and $Y$ as discussed earlier. Since
\linebreak $I\,\subseteq\,\ker p_{n(m)-1}^{'}\,=\,M_{n(m)}$, and
$I$ and $M_{n(m)}$ are generated by the same generators, we have
$I\,=\,M_{n(m)}$, that is, $\ker p_m\,=\,M_{n(m)}(A)$, where
$$M_{n(m)}(A)\,:=\,\{f\,\in\,A\,:\,o(f)\,\geq\,n(m)\}.$$ But then
$$\ker
q_m\,=\,M_{n(m)}(A_1)\,:=\,\{f\,\in\,A_1\,\subset\,A\,:\,o(f)\,\geq\,n(m)\},$$
where $A_1$ is as in Lemma 3.2. By Theorem 3.4, $A_1\,=\,\F_2$. It
follows that $A\,=\,\F_2$ topologically in view of the open mapping
theorem.

Case 2. for each $m$, $\ker p_m$ is singly generated by the
monomial $X^{n(m)}$. In this case
$$A_1\,=\,\{f\,\in\,\F_2\,:\,q_m(f)\,=\,\sum_{j=0}^{\infty}\sum_{i=0}^m\,\mid\lambda_{i,j}\mid\,<\,\infty\;\textrm{for
all}\;m\}\,=\,A_X,$$ by Theorem 3.4, since $X^{n(m)}\,\in\,\ker
q_m$. So, $A_X\,=\,A_1\,\subset\,A\,\subset\,\F_2$, by Lemma 3.2.
But, as we discussed in the property (5), $A_X\,\subset
\,A\,\subset\,A_X\,\subset\,\F_2$ since $p_m$ is a proper seminorm
on $A$ for each $m$ such that $X^{n(m)}\in \ker p_m$. Thus we have
$A\,=\,A_X$ topologically in view of the open mapping theorem.

Case 3. for each $m$, $\ker p_m$ is singly generated by the
monomial $Y^{n(m)}$. Follow the argument of Case 2, and the proof
is complete.

As corollaries, we have the following characterizations of a
Beurling-Fr\'{e}chet algebra $\ell^1(\zz^{+k},\,\Omega)$ of
semiweight type as a Fr\'{e}chet algebra.
\begin{cor} \label{Corollary 3.4_SM} Let $A$ be a Fr\'{e}chet algebra
of power series in $\F_k$. Then $A$ is equal to a
Beurling-Fr\'{e}chet algebra $\ell^1(\zz^{+k},\,\Omega)$ of
semiweight type if and only if the Fr\'{e}chet topology of $A$ is
defined by a sequence $(p_m)$ of proper seminorms. In particular,
$A\,=\,\F_k$ if and only if the Fr\'{e}chet topology of $A$ is
defined by a sequence $(p_m)$ of proper seminorms with
finite-dimensional cokernels.
\end{cor}

In fact, we have the following result on an Arens-Michael
representation of $A$.
\begin{cor} \label{Corollary 3.5_SM} Let $A$ be a Fr\'{e}chet algebra
of power series in $\F_k$ such that the polynomials are dense in
$A$. Then $A$ is not equal to a Beurling-Fr\'{e}chet algebra
$\ell^1(\zz^{+k},\,\Omega)$ of semiweight type if and only if
$A\,=\,\lim\limits_{\longleftarrow}A_{m}$ , where each $A_m$ is a
Banach algebra of power series in $\F_k$.
\end{cor}
{\it Proof.} Suppose that $A$ is not equal to a
Beurling-Fr\'{e}chet algebra \linebreak
$\ell^1(\zz^{+k},\,\Omega)$ of semiweight type. Evidently, by
Corollary \ref{Corollary 3.4_SM}, we may suppose that each $p_m$
is a norm on $A$. Now, by Theorem \ref{Theorem 3.1_SM}, $A$
contains a closed maximal ideal $M\,=\,\ker \pi_0$ such that
$$\bigcap_{n\geq 1}\overline{M^n}\,=\,\{0\}\; \textrm{and}
\;\textrm{dim}(\overline{M^n}/\overline{M^{n+1}})\,=\,C_{n+k-1,n}\;(n\,\in\,\nn).$$
By Proposition \ref{Proposition 2.3_SM}, for each sufficiently
large $l$, $M_l$ is a non-nilpotent maximal ideal of $A_l$ such
that
$$\bigcap_{n\geq 1}\overline{M_l^n}\,=\,\{0\}\;\textrm{and}\;
\textrm{dim}(\overline{M_l^n}/\overline{M_l^{n+1}})\,=\,C_{n+k-1,n}\;(n\,\in\,\nn).$$
Again, by Theorem \ref{Theorem 3.1_SM}, $A_l$ is a BAPS in $\F_k$
for each sufficiently large $l$. Hence, by passing to a suitable
subsequence of $(p_m)$ defining the same Fr\'{e}chet topology of
$A$, we conclude that each $A_l$ is a BAPS in $\F_k$.

The converse has already been discussed in the property (2) above.

The immediate consequence of Corollary ~\ref{Corollary 3.5_SM} is:
if $A$ is a FrAPS in $\F_k$ such that the polynomials are dense in
$A$ and such that it is not a Beurling-Fr\'{e}chet algebra
$\ell^1(\zz^{+k},\,\Omega)$ of semiweight type, then $A$ satisfies
Loy's condition (E) by Corollary 3.6. A somewhat more elaborate
version of the same idea enables us to drop the condition on the
polynomials in order to get a more general result, given below.

We recall that if the Fr\'{e}chet topology $\tau$ of $A$ is given
by a sequence $(p_m)$, then each $p_m$ is of type (E) if and only
if $A$ satisfies Loy's condition (E). Also, by a
several-variable-analogue of ~\cite[Theorem 2]{8}, $A$ satisfies
Loy's condition (E) if and only if $A$ admits a growth sequence;
i.e. there exists a sequence $(\sigma _K)_{K\,\in\,\nn^k}$ of
positive reals such that $\sigma _K\,\pi_K(x)\, \rightarrow \,0$
for each $x\, \in\, A$. This characterization of Loy's condition
(E) tells us that a Beurling-Fr\'{e}chet algebra
$\ell^1(\zz^{+k},\,\Omega)$ of semiweight type does not satisfy
Loy's condition (E) since $\F_r,\,1\,\leq\,r\,\leq\,k - 1$, does
not satisfy Loy's condition (E), which is a closed subalgebra by
the property (5).
\begin{thm} \label{Theorem 3.6_SM} Let $A$ be a Fr\'{e}chet algebra of power series in
$\F_k$. Then $A$ is not equal to the Beurling-Fr\'{e}chet algebra
$\ell^1(\zz^{+k},\,\Omega)$ of semiweight type if and only if
$A\,=\,\lim\limits_{\longleftarrow}A_{m}$ , where each $A_m$ is a
Banach algebra of power series in $\F_k$. In particular, $A$
satisfies Loy's condition (E) in this case.
\end{thm}
{\it Proof.} The proof is similar to that of ~\cite[Theorem
3.6]{12}, and so we will merely sketch it only. Supposing $A$ is
not equal to a Beurling-Fr\'{e}chet algebra
$\ell^1(\zz^{+k},\,\Omega)$ of semiweight type, we may suppose
that each $p_m$ is a norm on $A$, by Corollary ~\ref{Corollary
3.4_SM}. We first show that the projections $\pi_J$ are continuous
on $(A,\,p_m)$ for all $J$ and $m$; i.e. the inclusion maps
$(A,\,p_m) \hookrightarrow \F_k$ is continuous.

Since $M \,=\, \ker \pi_{0,\,\dots,\,0}$ is a non-nilpotent,
closed maximal ideal of $A$ such that $\bigcap _{n\geq 1}
\overline{M^n}\,=\,\{0\}$, $M$ is a non-nilpotent, maximal ideal
of $(A,\,p_m)$ for each $m$; also,
$\overline{M^{n+1}}\,\neq\,\overline{M^{n}}\,\neq\,\{0\}$ (closure
in $(A,\,\tau)$) in $(A,\,p_m)$ for all $m$ and $n$. In fact, we
may suppose that for each sufficiently large $m$, $M$ is also
closed in $(A,\,p_m)$ and hence that $\pi_{0,\,\dots,\,0}$ is
$p_m$-continuous for each sufficiently large $m$. Also, by the
argument in the proof of ~\cite[Proposition 2.3]{12}, we may
suppose that $M_m$ is a non-nilpotent, maximal ideal of $A_m$ and
that $\bigcap_{n\geq 1}\overline{M_m^n}\, =\, \{0\}$ for each $m$.
Assume inductively that $\pi_I$ is continuous for each $\mid I
\mid\,<\,\mid J \mid$, and take $(x_n)$ in $(A,\,p_m)$ with
$p_m(x_n) \,\rightarrow \,0$. Then, following the argument given
in Theorem ~\ref{Theorem 3.1_SM}, we deduce that some non-zero
linear combination of $X^I,\,\mid I\mid\,=\,\mid J\mid$, lies in
$A \bigcap \overline{M_m^{\mid J\mid+1}}$ (which is, in fact,
$\overline{M^{\mid J\mid+1}}$ in $(A,\,p_m)$ for each $m$), a
contradiction of the fact that $\overline{M_m^{\mid J\mid+1}}
\,\neq\, \overline{M_m^{\mid J\mid}}$.

Next, for each $m\,\in\,\nn$, one shows that $p_m$ is closable on
$A$ by noticing that the inclusion map $(A,\,p_m) \,\hookrightarrow
\,\F_k$ can be extended to a continuous homomorphism $\phi : A_m
\rightarrow \,\F_k$. So $x_n \rightarrow x$ in $(\F_k,\,\tau_c)$,
and hence $p_m(x_n) \rightarrow 0$.

Then one shows that $\phi$ is, indeed, injective, and so $A_m$ is a
Banach algebra of power series in $\F_k$ for each $m$. In
particular, $A$ satisfies Loy's condition (E), by the remark
following Corollary ~\ref{Corollary 3.2_BP}.
\section{Automatic continuity and uniqueness of\\ topology, and open questions.} We
now turn to establish that every FrAPS in $\F_k$ has a unique
Fr\'{e}chet topology. By ~\cite[Theorem 1]{9}, it is clear that
every FrAPS $A$ in $\F_k$ satisfying Loy's condition (E) has a
unique Fr\'{e}chet topology. Since a Beurling-Fr\'{e}chet algebra
$\ell^1(\zz^{+k},\,\Omega)$ of semiweight type does not satisfy
the condition (E), Theorem ~\ref{Theorem 3.6_SM} gives the
following result from ~\cite{10} by noticing that ~\cite[Theorem
10]{10} holds true if $A$ is considered to be a FrAPS in $\F_k$ in
the codomain.
\begin{thm} \label{Theorem 4.1_SM} Let $A$ be a Fr\'{e}chet algebra of power series
in $\F_k$ such that $A$ is not equal to a Beurling-Fr\'{e}chet
algebra $\ell^1(\zz^{+k},\,\Omega)$ of semiweight type. Then a
homomorphism $\theta \,:\, B\,\rightarrow\,A$ from a Fr\'{e}chet
algebra $B$ is continuous provided that the range of $\theta$ is
not one-dimensional.
\end{thm}

Now, we prove the uniqueness of the Fr\'{e}chet topology for
Beurling-Fr\'{e}chet algebras $\ell^1(\zz^{+k},\,\Omega)$ of
semiweight type.
\begin{thm} \label{Theorem 4.2_new} Let $A$ be a Beurling-Fr\'{e}chet
algebra $\ell^1(\zz^{+k},\,\Omega)$ of semiweight type. Then $A$
has a unique Fr\'{e}chet topology.
\end{thm}
{\it Proof.} Let $\tau^{'}$ be another topology such that $(A,
\tau^{'})$ is a Fr\'{e}chet algebra. We recall that the original
Fr\'{e}chet topology $\tau$ of $A$ is defined by a sequence
$(p_m)$ of proper power series seminorms. First we note that $A$
is a local Fr\'{e}chet algebra, and hence it has a unique maximal
ideal $M = \ker \pi_{0,\dots,\,0}$ (this is clearly an algebraic
property). Since $(A,\,\tau{'})$ is local, $A$ is a $Q$-algebra by
~\cite[Corollary 3]{18}. Hence $M$ is $\tau^{'}$-closed. For
$n\,\in\,\zz^+$, let
$$M_n(A)\,:=\,\{f\,\in\,A\,:\,o(f)\,\geq\,n\}.$$ Since $(A,\,\tau)$ is
a FrAPS in $\F_k$ such that the polynomials are dense in $A$, by
argument of Theorem 3.1, we have $$\bigcap_{n\geq
1}\overline{M^n}^\tau = \{0\}\;\textrm{and} \;M_n(A) =
\overline{M^n}^\tau,$$ for each $n$, hence $\bigcap_{n\geq
1}M_n(A) = \{0\}$. Now let $B$ be the closure of the polynomials
in $k$ variables in the topology $\tau{'}$. Hence $B$ is a closed
subalgebra of $(A, \tau^{'})$. Then $M \bigcap B$ is a
non-nilpotent, closed maximal ideal of $B$ such that
$\bigcap_{n\geq 1}\overline{(M\,\bigcap\,B)^n}^{\tau^{'}} =
\{0\}$; the latter holds because $$\bigcap_{n\geq 1}M_n(B) =
\{0\}\;\textrm{and}\;(M\,\bigcap\,B)^n \subset
\overline{(M\,\bigcap\,B)^n}^{\tau^{'}} \subset M_n(B),$$ for each
$n$. Also, since the polynomials in $k$ variables are dense in
$B$, $\overline{(M\,\bigcap\,B)^n}^{\tau^{'}} = M_n(B)$, and so
$$\textrm{dim}(\overline{(M\,\bigcap\,B)^n}^{\tau^{'}}/\overline{(M\,\bigcap\,B)^{n+1}}^{\tau^{'}})
= C_{n+k-1,n},$$ for all $n$. Now, by Theorem 3.1, $(B,\,\tau^{'})$
is a FrAPS in $\F_k$.

Since the inclusion map $(B,\,\tau{'})\,\hookrightarrow
\,(A,\,\tau{'})$ is continuous, by the open mapping theorem for
Fr\'{e}chet spaces, the inclusion map is a linear homeomorphism.
Hence both $(A,\,\tau)$ and $(A,\,\tau{'})$ are FrAPS in $\F_k$.
Let $(a_n)$ be a sequence in $A$ such that $a_n\,\rightarrow \,0$
in $(A,\,\tau)$ and $a_n\,\rightarrow \,a$ in $(A,\,\tau{'})$. For
each $J\,\in\,\nn^k$, the functional $\pi_J$ is continuous on both
$(A,\,\tau)$ and $(A,\,\tau{'})$, and so $\pi_J(a)\,=\,0$, whence
$a\,=\,0$. By the closed graph theorem for Fr\'{e}chet spaces, the
identity map $(A,\,\tau)\,\rightarrow\,(A,\,\tau{'})$ is a linear
homeomorphism, and so $\tau\,=\,\tau{'}$ on $A$. The proof is
complete.

As a corollary of the above two results, we have the following
result (and Questions) in the theory of automatic continuity. We
recall that the continuity of automorphism of $\F_k$ is proved in
~\cite[p. 136]{17}, and the uniqueness of the Fr\'{e}chet topology
of $\F_k$ is proved in ~\cite[Theorem 4.6.1]{3}. The following
gives an answer to ~\cite[Question 11]{2} for FrAPS in $\F_k$.
\begin{cor} \label{Corollary 4.2_SM} Let $A$ be a Fr\'{e}chet algebra of power series
in $\F_k$. Then $A$ has a unique Fr\'{e}chet topology.
\end{cor}

We now list the following questions, which may have some interest in
the theory of automatic continuity.

\noindent {\bf Question 1.} Let $A$ be a Beurling-Fr\'{e}chet
algebra $\ell^1(\zz^{+k},\,\Omega)\,(\neq\,\F_k)$ of semiweight
type. Is every automorphism of $A$ continuous?

The author conjecture that the answer of Question 1 is in the
affirmative; also, this question has been studied for special cases
(see ~\cite{5, 6}). More generally, we have the following

\noindent {\bf Question 2.} Is every homomorphism
$\theta\,:\,B\,\rightarrow\,\F_k\;(k\,>\,1)$ from a Fr\'{e}chet
algebra $B$ continuous?

The above question has recently been settled partially in ~\cite{4}:
(i) we remark that if $\theta\,:\,B\,\rightarrow\,\F_k$ is a
discontinuous homomorphism from an $(F)$-algebra $B$ (i.e., a
complete, metrizable topological algebra) into $\F_k$ such that
$\theta(B)$ is dense in $(\F_k,\,\tau_c)$ and such that the
separating ideal of $\theta$ has finite codimension in $\F_k$, then
$\theta$ is an epimorphism (see ~\cite[Theorem 12.1]{4}), and (ii)
there exists a Banach algebra $(A,\,\|\cdot\|)$ such that
$\cc[X_1,\,X_2]\,\subset\,A\,\subset\,\F_2$, but such that the
embedding $(A,\,\|\cdot\|)\,\hookrightarrow\,(\F_2,\,\tau_c)$ is not
continuous (see ~\cite[Theorem 12.3]{4}); surprisingly, $A$ is
(isometrically isomorphic to) a Banach algebra of power series in
$\F_1$ (see ~\cite[Theorem 10.1 (i)]{4}).

As shown in ~\cite{14}, $\F_\infty$ admits two inequivalent
Fr\'{e}chet algebra topologies. It is clear that the unique
maximal ideal $M\,=\,\ker\,\pi_0\; (0\,=\,(0,\,0,\,\dots)),$ is
closed under both topologies, since $\F_\infty$ is a local
Fr\'{e}chet $Q$-algebra. Thus it is of interest to know whether
every FrAPS in $\F_\infty$ (except $\F_\infty$ itself) has a
unique Fr\'{e}chet topology. In other words, we have the following
natural question.

{\bf Question 3.} Is there any other proper, unital subalgebra of
$\F_\infty$, with two inequivalent Fr\'{e}chet topologies? In
particular, is there any other proper subalgebra of $\F_\infty$
which is closed under the topology imposed by Charles Read on
$\F_\infty$ and which is also FrAPS in $\F_\infty$ in its ``usual"
topology?

To answer the latter part of the above question, the ``natural"
extension of Beurling-Fr\'{e}chet algebras of semiweight type
(i.e., $\ell^1((\zz^+)^{<\,\omega},\,\Omega)$) would be an easy
target. Also, we have the following curious question.

{\bf Question 4.} Does there exist a Fr\'{e}chet algebra with
infinitely many inequivalent Fr\'{e}chet algebra topologies?

For $m\,\in\,\nn$, set $${\mathcal
U}_m=\{f=\sum\{\alpha_rX^r:r\in(\zz^+)^{<\omega}\}\in\F_\infty:p_m(f):=\sum\mid\alpha_r\mid
\,m^{\mid r\mid}<\infty\},$$ and then set $${\mathcal
U}\,=\,\bigcap\{{\mathcal U}_m\,:\,m\,\in\,\nn\}.$$ It is clear
that each $({\mathcal U}_m,\,p_m)$ is a unital Banach subalgebra
of $\F_\infty$ and that $({\mathcal U},\,(p_m))$ is a unital
Fr\'{e}chet subalgebra of $\F_\infty$. Being semisimple
Fr\'{e}chet algebras, the test algebra $\cal U$ for the still
unsolved ``Michael problem" (and ${\cal U}_{m}$ for each $m$
appearing in the Arens-Michael representation of $\cal U$),
$\ell^1((\zz^+)^{<\omega})$ and the algebra $\ell^1(S_c)$, where
$S_c$ denotes the free semigroup on $c$ generators, have unique
Fr\'{e}chet topologies (the result also follows from
~\cite[Theorems 10.1 and 10.5]{4} and ~\cite[Corollary 4.2]{12}).
In this regard, we first give the following result (for the proof,
see either ~\cite[Theorem 1]{8} or ~\cite[Theorem 2]{9}) by
noticing that FrAPS $A$ in $\F_\infty$ satisfies Loy's condition
(E) if there is a sequence
$(\gamma_K\,:\,K\,\in\,\nn^k,\,k\,\in\,\nn)$ of positive reals
such that $(\gamma_K^{-1}\pi_K)$ is equicontinuous.
\begin{thm} \label{Theorem_Loy} Let $A$ be a Fr\'{e}chet algebra of
power series in $\F_\infty$ satisfying Loy's condition
(E), above, and let $\phi \,:\,B\,\rightarrow \,A$ be a
homomorphism from a Fr\'{e}chet algebra $B$ into $A$ such that
$X_1\,\in\,\phi(B)$. Then $\phi$ is continuous. In particular,
every automorphism of $A$ is continuous, and $A$ has a unique
Fr\'{e}chet algebra topology.
\end{thm}

We remark that the Fr\'{e}chet topology of $A$ of the above theorem
is defined by a sequence of norms. We have the following ``natural"
generalization of ~\cite[Theorem 10.1]{4}.
\begin{thm} \label{Theorem_10.1} Let $A$ be a Fr\'{e}chet algebra of power series in $\F_\infty$,
with its topology defined by a sequence $(p_m)$ of norms. Suppose
that $A$ is a graded subalgebra of $\F_\infty$. Then there is a
continuous embedding $\theta$ of $A$ into $(\F, \,\tau_c)$ such that
$\theta(X_1)\,=\,X$, and so $A$ is (isometrically isomorphic to) a
Fr\'{e}chet algebra of power series in $\F_1$. In particular, $A$
has a unique Fr\'{e}chet topology.
\end{thm}

We note that we cannot drop the assumption on $(p_m)$ in the above
theorem; $(\F_\infty, \,\tau_c)$ is a counterexample, since it can
continuously be embedded in $(\F_2, \,\tau_c)$, by ~\cite[Theorem
9.1]{4}, but, by Theorem 2.6 (or Theorem 11.8) of ~\cite{4}, there
is no embedding of $\F_2$ into $\F$.

Finally, we remark that we have recently established the uniqueness
of (F)-algebra topology for (F)-algebra of power series in the
indeterminate $X$ (see ~\cite[Corollary 11.7]{4}). Then the natural
question is to extend this result for the several-variable case, in
order to settle the Dales question completely. But our approach
fails, because the clause of the continuity of coordinate
projections cannot be dropped from the definitions of Banach and
Fr\'{e}chet algebras of power series in $k$ indeterminates (see
~\cite[Theorem 12.3]{4}). The approach here is based on a sequence
$(p_m)$ of submultiplicative seminorms, and so these results cannot
be extended to {\it all} (F)-algebra of power series in $k$
indeterminates.

{\bf Acknowledgements.} The author would like to thank Professors
H. G. Dales for encouraging me to solve this problem and C. J.
Read for bringing $A_X$ (which does not satisfy Loy's condition
(E)) to my attention. This work was partially supported by the
Commonwealth Scholarship Commission [INCF-2008-63].

\end{document}